\documentclass{article}
\usepackage[utf8]{inputenc}
\usepackage{listings}
\usepackage{booktabs}
\usepackage{tabularx}
\usepackage{multirow}
\usepackage{amsmath}
\usepackage{amssymb}
\usepackage{longtable}
\usepackage{url}
\usepackage{cite}
\usepackage{graphicx}
\usepackage{ulem}
\usepackage{floatrow}
\usepackage{caption}
\usepackage{authblk}
\newfloatcommand{capbtabbox}{table}[][\FBwidth]
\usepackage{xcolor}
\usepackage{multicol}
\usepackage{setspace}
\usepackage{gensymb}
\usepackage{geometry}
\usepackage{graphicx}
\usepackage{multirow}
\usepackage{algorithm}
\usepackage{algpseudocode}
\usepackage{multirow}
\usepackage{csquotes}
\usepackage{amsmath}
\usepackage{bm}
\usepackage{natbib}  

\usepackage{adjustbox}
\usepackage{titlesec}

\usepackage{tabularx,ragged2e,booktabs,caption}
\newcolumntype{C}[1]{>{\Centering}m{#1}}

\newcommand{\RN}[1]{%
	\textup{\uppercase\expandafter{\romannumeral#1}}%
}
\usepackage{tabularx,ragged2e,booktabs,caption}
\usepackage{colortbl,hhline}
\definecolor{c1}{rgb}{0.30980, 0.50588, 0.73725}
\definecolor{c2}{rgb}{0.82353, 0.87843, 0.92941}
\usepackage{adjustbox}
\usepackage{textcomp} 
\usepackage{subfig}
\titleformat{\paragraph}
{\normalfont\normalsize\bfseries}{\theparagraph}{1em}{}
\titlespacing*{\paragraph}
{0pt}{3.25ex plus 1ex minus .2ex}{1.5ex plus .2ex}
\usepackage{multirow}

\usepackage{ragged2e}
\renewcommand{\raggedright}{\leftskip=0pt \rightskip=pt plus 0cm}

\usepackage{amssymb}
\usepackage{gensymb}
\usepackage[toc,page]{appendix}

\title{Isogeometric contact analysis in subsea umbilical and power cables}

\author[a,b]{Tianjiao Dai}
\author[c]{Shuo Yang}
\author[d]{Xing Jin \footnote{corresponding author, Email: Jinxing@coffshore.cn}}
\author[e]{Svein Sævik} 
\author[a]{Jiaxuan Zhang}
\author[a]{Jun Wu}
\author[f]{Naiquan Ye}
\affil[a]{School of Naval Architecture and Ocean Engineering, Huazhong University of Science and Technology (HUST), 430074, Wuhan, China}
\affil[b]{Hubei Key Laboratory of Naval Architecture \& Ocean Engineering Hydrodynamics,HUST, 430074, Wuhan, China}
\affil[c]{Zoomlion Heavy Industry Science and Technology Co., Ltd., China}
\affil[d]{China Offshore Engineering \& Technology Co., Ltd., China}
\affil[e]{Department of Marine Technology, Norwegian University of Science and Technology, NO-7491 Trondheim, Norway}
\affil[f]{Energy and Transport, Sintef Ocean, NO-7052 Trondheim, Norway}

\begin{document}
\date{}
\maketitle
\setstretch{1.5}

\section*{Abstract}

Subsea umbilical and power cables contain a large number of contact interfaces between different geometries and materials. These complex interactions rise significant challenges for accurately considering contact surface properties by using traditional analytical solutions or finite element methods. These properties have been identified as the most sensitive parameters when performing the numerical simulation for stress analysis. Therefore, it is essential to apply a novel approach for contact analysis which improves the accuracy and efficiency for predicting contact properties. This paper presents an isogeometric analysis (IGA) approach addressing contact problems in dynamic umbilicals and power cables. Firstly, this isogeometric contact algorithm is formulated in MATLAB as a tool including the geometry description, contact detection and penalty function. Secondly, the contact interface between a steel tube and an outer sheath in an dynamic umbilical is established by this IGA contact algorithm and validated against that in ABAQUS for proving the accuracy and efficiency of IGA. Finally, the effects of element refinement, geometrical description, penalty factor on the accuracy, efficiency and stability of IGA are discussed. \\ 

Keywords: isogeometric contact, dynamic umbilical, contact force, contact stiffness

\section{Introduction}
\noindent

The subsea umbilical cable is a highly integrated system serving as lifelines for supplying power and controls to subsea exploration and production equipments. It is one of the most important and vulnerable parts in the subsea production system when especially moving to the ultra-deep water. They burden the alternating tension and bending moments along the whole length of the cable resulting in fatigue damage accumulated in the helical components. The fatigue damage is calculated by using the stress range of helical components which is determined by the contact force along the helical components. It is well known that this contact force is affected by the contact surface's hardness since it would directly determine the radial motion and further influences the tensile capacity of the corresponding helical component.  Therefore, it is crucial to propose a reliable methodology for contact analysis in subsea umbilical and power cables. \\

\begin{figure}[htbp]
\centering
\includegraphics[width=0.8\linewidth]{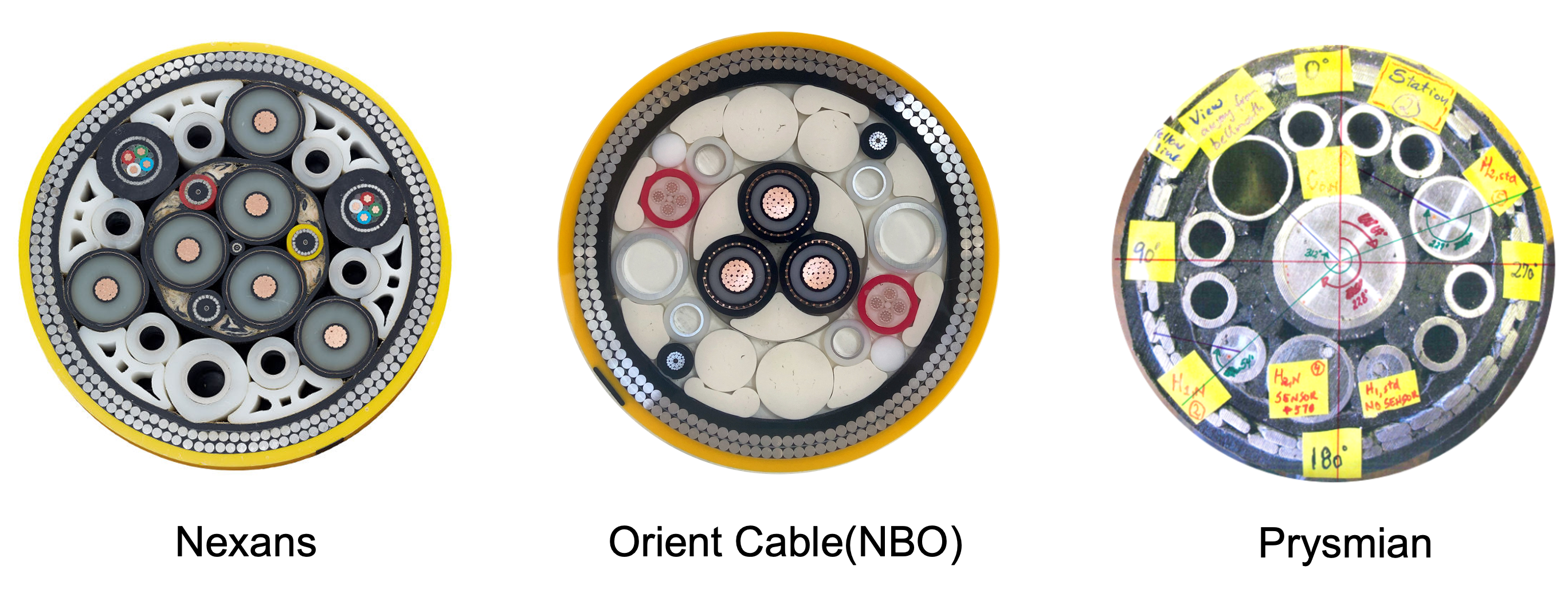}
\caption{ Three typical cross sections of dynamic umbilical cables and these umbilicals are from Nexans, Orient Cable(NBO) and Prysmian, respectively.}
\label{fig:umbi-cross-section}
\end{figure}
Typical cross sections of dynamic umbilical cables are presented in Figure \ref{fig:umbi-cross-section}. The internal contacts occur between various components including fillers, cables, copper wires, sheath and steel tubes and result in contact pairs between different materials and geometries. These contacts rise significant challenges in modelling and convergence for the numerical analysis. Until now, three main methods are applied to perform the stress analysis of dynamic umbilicals. The main differences among these methods lie in the approaches how to simplify the geometry of helical components and how to model the contact mechanism in the cross section. The first method is to establish finite element models in ABAQUS or ANSYS and is the most widely applied method in the engineering industry. The contact is primarily modelled by using the penalty function \cite{li_efficient_2022,li_finite_2021,li_finite_2023,wang_full_2021,wang_novel_2022,zhang_research_2020,yin_experimental_2021,yang_study_2021,qin_review_2024}. However, this factor is typically chosen either by default or obtained through a sensitivity study according to the numerical convergence performance. This factor does not accurately describe the real contact surface's hardness since it fails to reflect these contact pairs' material properties and geometry effects on the contact hardness. If the solid element is applied together with this penalty factor, the contact properties could be caught. In addition, such numerical analysis costs enormous computation time which is not preferred in industry during the design phase. The second method contains some simplified analytical models \cite{saevik_theoretical_2011,witz_axial-torsional_1992} where no gap is assumed between neighbouring layers. This enforces the radial deformation coordination conditions when solving the equation system. This means the radial motion is only determined by the thickness change of each layer rather than the contact surface properties. For the third method, the Sandwich beam theory is employed and curved beam is used to model the internal helical components\cite{saevik_theoretical_2011}. This method enables efficient stress analysis of helical components due to a small number of element and the contact stiffness is used to describe the contact hardness\cite{saevik_stresses_1992}. The contact stiffness could be obtained from an analytical solution \cite{dai_anisotropic_2018,dai_friction_2017,dai_experimental_2020}. In their work, the axial stress behavior of the helical steel tube in the dynamic umbilical under 200KN tension closely matched experimental test results. However, deviations and convergence problems were found for low or high tension. This may be because the contact stiffness obtained by the analytical solution is not sufficient to reflect the geometry and material properties of contact interfaces. Therefore, it is important to propose an individual contact algorithm to calculate the real contact stiffness accounting for these nonlinear effects.\\

Until now, two primary methods are applied to establish the individual contact models. The first is Hertzian contact theory which is only valid when the contact area is sufficiently small\cite{wu_hertzian_2016,zhang_spherical_2014,guo_modified_2020}. Due to significant differences in Young's modulus between the sheath layer and the steel pipe, substantial axial loads leads to an increased contact area, as illustrated in Figure\ref{fig:umbi-cross-section}. This leads to the assumption of small contact areas in the Hertz contact theory no longer applicable, thus unsuitable for calculating contact stiffness within the umbilical cable. The second is the traditional finite element modeling in ABAQUS or ANSYS. In this approach, a series of straight lines are applied to approximate the curves describing the tube and sheath outlines. This is challenging to precisely capture the geometry details of the steel tube and sheath layer\cite{el-abbasi_modelling_2001}. If using the penalty function approach in the contact model, the actual penetration is smaller than that used in the numerical analysis. The penetration in the numerical analysis is calculated based on the straight line between nodes rather than the real curve, which leads to underestimated contact force. Additionally, the contact geometry outline at the nodes is not continuous resulting in different normal vector directions and penetrations from the adjacent elements. Consequently, the corresponding contact forces at the same node are discontinuous resulting in convergence problems\cite{de_lorenzis_mortar_2012,lu_circular_2009}. The accuracy and stability of numerical results could be improved by increasing the number of elements, however, this significantly increases computation times. Moreover, the prediction of contact surface's properties is highly sensitive to the the element number used. Hegadekatte et al.\cite{hegadekatte_multi_2010} performed contact analysis for planetary gear train to investigate wear behavior by using traditional finite element methods. Although the numerical results matched well with experimental data, the high-quality meshes are required for accurately describing contact behavior within planetary gear train and cost a significant amount of time for mesh generation. Additionally, the entire analysis process include iterations and re-meshing is required in each iteration. This was because the contact surfaces in planetary gear systems were complex and the contact boundaries could only be approximated by increasing the number of elements by using traditional finite element methods. These difficulties are unavoidable drawbacks of traditional finite element methods.\\

Zhao and Li\cite{zhao_solution_2011} established a finite element model for solving rolling friction contact problems. The numerical results align well with Hertzian theory with respect to the normal stress. However, due to the discretization inherent in traditional finite element methods, the contact surface radius can only be an integer multiple of the element size resulting in discrepancies with the Hertzian solution. Additionally, due to this discretization difficulty, sufficiently small elements had to be used for achieving accurate results resulting in long computational time. Singh et al.\cite{singh_analyses_2014} conducted a contact force analysis for a roller rolling over the raceways of a bearing. The finite element results showed that the overall trend of the contact force was almost consistent with the theoretical solution. However, the oscillations of contact force over time were observed. This was due to the polygonal approximation of the arc boundary. As the polygonal boundary rolls over the raceway, small impacts are generated and then the contact force oscillates. Some similar behaviors were found in Xing et al.\cite{xing_node--node_2019}'s work and Otto et al.\cite{otto_coupling_2019}'s study. To address these issues, they proposed to add artificial NURBS layers between two contact bodies to enhance the continuity of the contact surface boundary. This method significantly improved the results' accuracy without spurious oscillations. Therefore, this paper employed Isogeometric Analysis (IGA) to establish the contact interaction surface between the steel tube and the sheath layer.\\

This paper focuses on addressing 2D contact interaction surfaces in dynamic umbilical and power cables. The contact algorithms by IGA is introduced and formulated in MATLAB starting from the geometry basic functions, normal contact description, KTS (knot to surface) algorithm, contact detection and penalty function. Then, the same contact interaction surface is also analysed in ABAQUS to validate the accuracy, efficiency and stability of IGA contact algorithm. Moreover, the advantages of IGA is further discussed in terms of the element refinement, penalty factor and contact penetration.\\

\section{IGA contact algorithm}
\subsection{Typical contact pairs in dynamic umbilical}

\begin{figure}[htbp]
\centering
\includegraphics[width=0.6\linewidth]{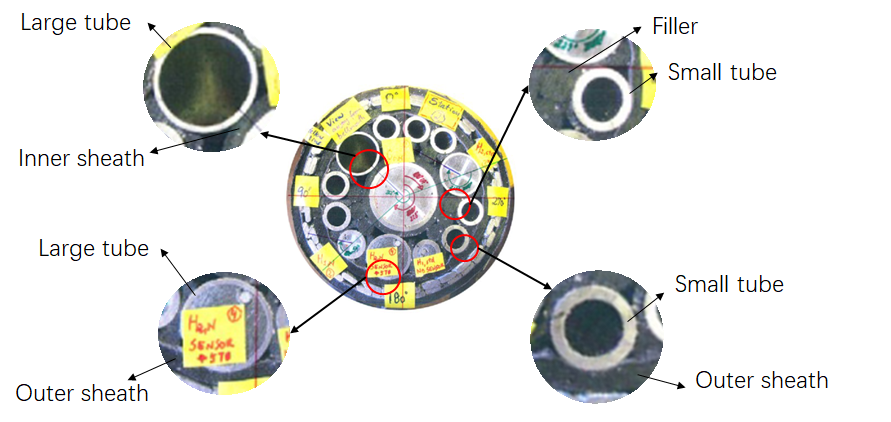}
\caption{The typical contact pairs in dynamic umbilicals}
\label{fig:typical contact pairs}
\end{figure}

Some typical contact pairs are selected to evaluate the performance of IGA contact algorithm with respect to its geometry description, accuracy, efficiency and convergence stability. Figure \ref{fig:typical contact pairs} presents these typical contact pairs in an umbilical specimen of a joint industry project conducted at SINTEF OCEAN. A detailed introduction of the full-scale test can be found in \cite{dai_experimental_2020}. These contact pairs are also chosen for further verification with respect to the dynamic stress behavior as Part B of this study. In this study, the contact surface between the steel tube and sheath is examined to illustrate the IGA contact algorithm's capacity in geometry description and multi-material considerations. The following subsections will present how these contact geometries are described by using NURBS and how the contact stiffness is determined by using the penalty method. \\

\subsection{Geometry description}
The contact pair of outer sheath and large tube is taken as an example to illustrate the isogeometric description by using NURBS curve. The B-spline curve is firstly introduced in order to understand the NURBS curve. 

\subsubsection{B-spline curve}
The B-spline curve is a piece-wise polynomial consisted by a linear combination of basis functions. Various methods could be applied to define these basis functions, but in this study, the most commonly used recurrence formula will be presented. Define a knot vector $\mathbf{U} = \{u_1, u_2, \ldots, u_m\}$, where $u_i$ is denoted as a knot, and $\mathbf{U}$ constitutes a set of $m$ non-decreasing real numbers, i.e., $u_i \leq u_{i+1}$ for $i = 1, 2, \ldots, m-1$. The $N_{i,p}(u)$ symbolizes the $i$-th $p$-degree (order $P+1$) B-spline basis function. The specific recursive definition is expressed by Equation \ref{eq: eq1}:

\begin{equation}
\centering
\left\{\begin{array}{l}
N_{i, 0}(u)=\left\{\begin{array}{l}
1, \text { if } u_i \leq u<u_{i+1} \\
0, \text { else }
\end{array}\right. \\
N_{i, p}(u)=\frac{u-u_i}{u_{i+p}-u_j} N_{i, p-1}(u)+\frac{u_{i+p+1}-u}{u_{i+p+1}-u_{i+1}} N_{i+1, p-1}(u)
\end{array}\right.
\label{eq: eq1}
\end{equation}

When performing differentiation, the numerator and denominator may both be zero, then we define $0/0 = 0$. Additionally, it should be noted that the order $p$ of the basis function must not be exceeded by the differentiation order $k$; otherwise, the derivatives beyond order $p$ are all zero.\\

The B-spline curve consists of a linear combination of basis functions and control points. The p-th degree B-spline, controlled by n control points, is represented as follows:
\begin{equation}
\centering
 C(u) = \sum_{i=0}^{n} N_{i,p}(u)P_i
\label{eq: eq6}    
\end{equation}
Where, $N_{i,p}(u)$ is the $p$-th degree basis function defined on the non-uniform, non-periodic knot vector:
\[\mathbf{U} = \{\underbrace{a, \ldots, a}_{p+1}, u_{p+2}, \ldots, u_{m-p-1}, \underbrace{b, \ldots, b}_{p+1}\}\]
$P_i$ represents the control points. Typically, $a=0$ and $b=1$. The details of important properties of the B-spline basis functions is found in \cite{lu_isogeometric_2011}. In summary, the B-spline curves can provide strong local support and controllable continuities. This makes them flexible and accurate to describe the specified curve geometry. \\

\subsubsection{NURBS curve}
NURBS curves are rational B-splines and capable of flexibly and accurately representing various curves. This capacity can not be achieved by polynomial curves. NURBS curves incorporate additional weights, denoted by $w$, to describe the relative influence of the control points. A p-th degree NURBS curve is defined by: 
\begin{equation}
C(u)=\sum_{i=1}^n R_{i, p}(u) P_i, \quad a \leq u \leq b
\label{eq:NURBS}
\end{equation}

\begin{equation}
R_{i, p}(u)=\frac{N_{i, p}(u) w_i}{W(u)}=\frac{N_{i, p}(u) w_i}{\sum_{j=0}^n N_{i, p}(u) w_j}
\end{equation}

where $N_i,p(u)$ are the p-th degree B-spline basis functions, $P_i$ are the control points, and $w_i$ are the corresponding weights. \\

\begin{figure}[htbp]
\centering
\includegraphics[width=0.6\linewidth]{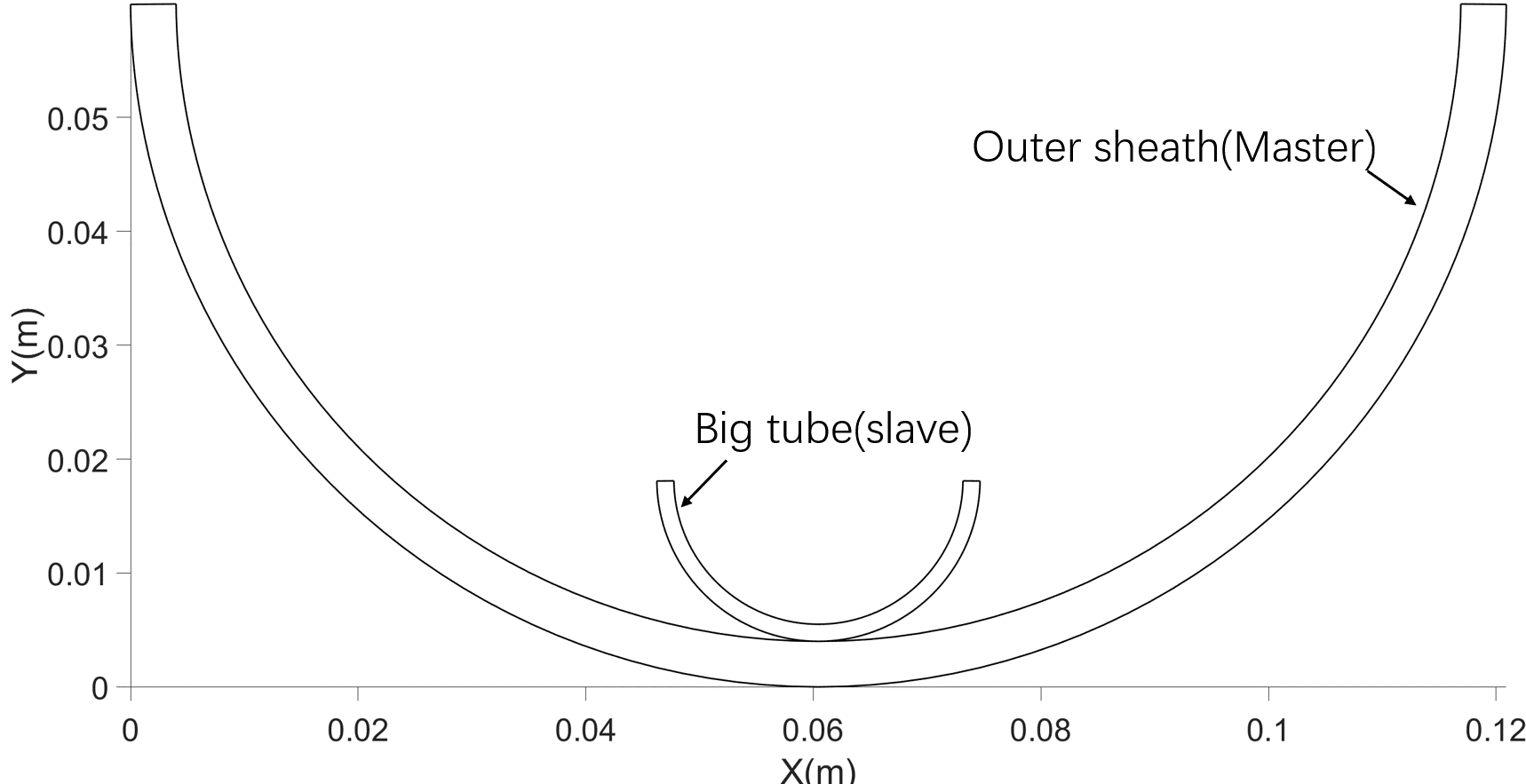}
\caption{The contact pair geometry of outer sheath and large tube}
\label{fig:actualgeometry}
\end{figure}

\begin{figure}[htbp]
    \centering
    \includegraphics[width=0.8\linewidth]{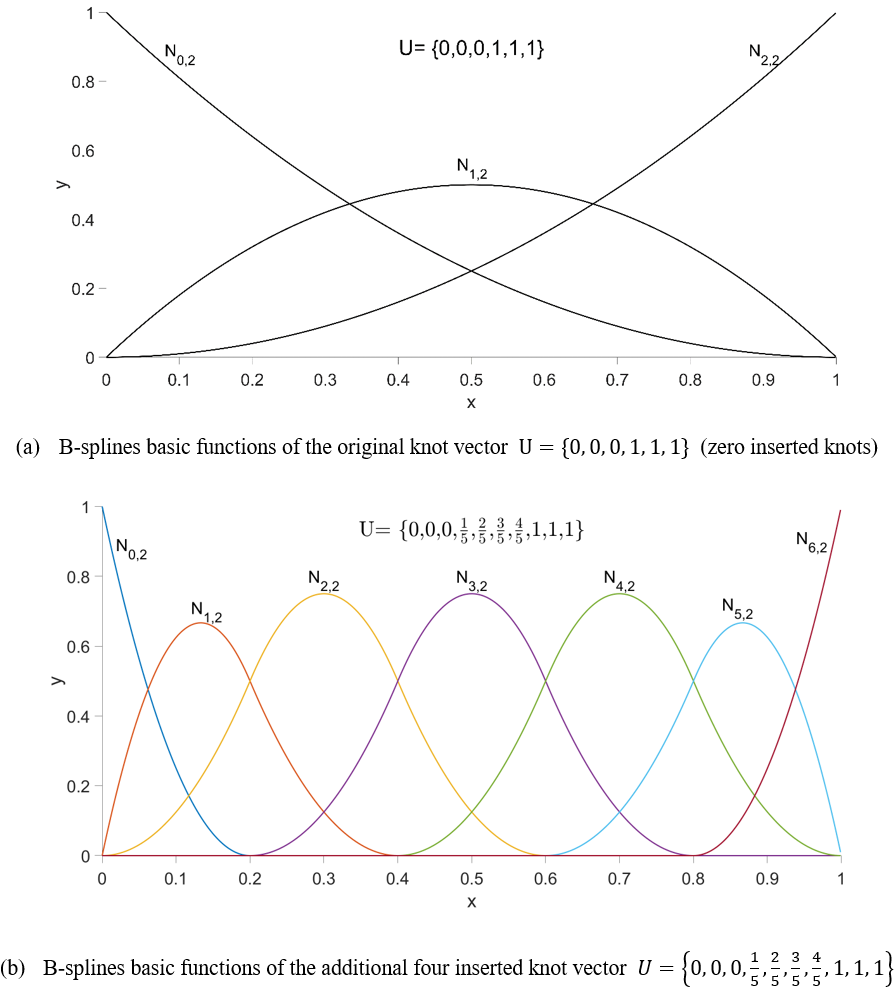}
    \caption{B-splines basic functions of zero inserted knots and four inserted knots}
    \label{fig:0and4knotsfunctions}
\end{figure}

\begin{figure}[htbp]
    \centering
    \includegraphics[width=\linewidth]{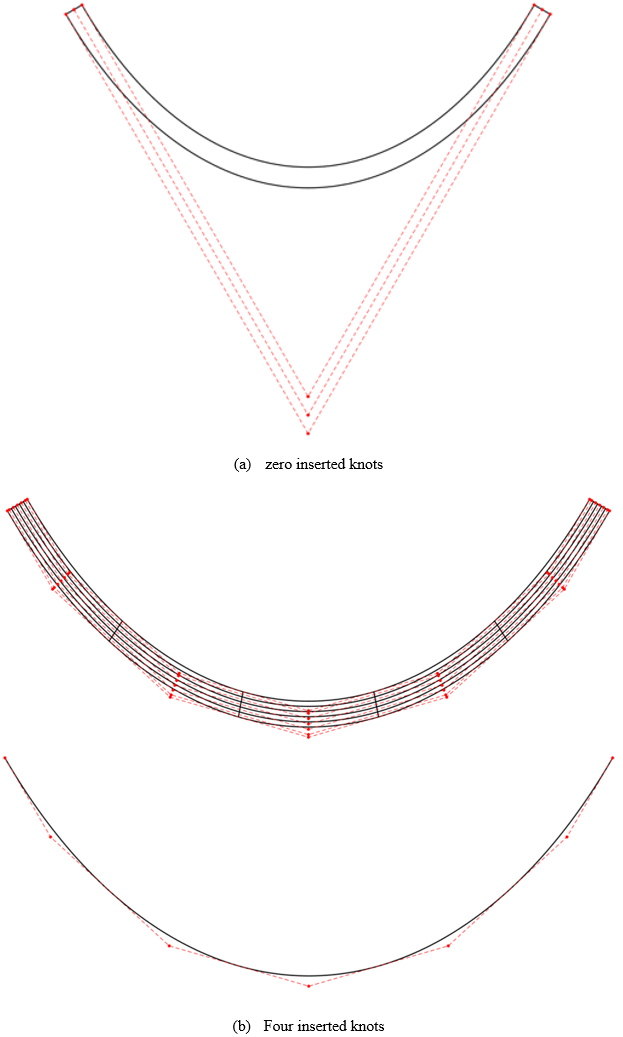}
    \caption{A 120\degree  arc of the outer sheath described by using zero inserted knots and four inserted knots}
    \label{fig:0and4knotsmesh}
\end{figure}

Figure \ref{fig:actualgeometry} illustrates the real geometries of outer sheath and steel tube as introduced above. Two different sets of knots are used in order to evaluate the geometry description capacity. The NURBS curve is constructed to span from the start to the end nodes of the knot vector. The original knot vector is given by $\mathbf{U}=\{0,0,0,1,1,1\}$. Figure \ref{fig:0and4knotsfunctions} shows the basic functions over the original knot vector and those over the knot vector with four inserted knots, which are generated according to B-splines. The four inserted knots are placed between $\{u_2,u_3\}$ and a new knot vector is formed as $\mathbf{U}=\{0,0,0,\frac{1}{5},\frac{2}{5},\frac{3}{5},\frac{4}{5},1,1,1\}$.  A 120\degree arc is selected in order to visibly illustrate the effect of control points. The corresponding generated geometry is presented in Figure \ref{fig:0and4knotsmesh}. It is evident that the circle can be accurately described even by the original knot vector. This means the number of elements used in the finite element method will be significantly reduced by using the NURBS curve to describe the contact boundaries. \\

\subsection{Normal contact}
\begin{figure}[htbp]
\centering
\includegraphics[width=0.6\linewidth]{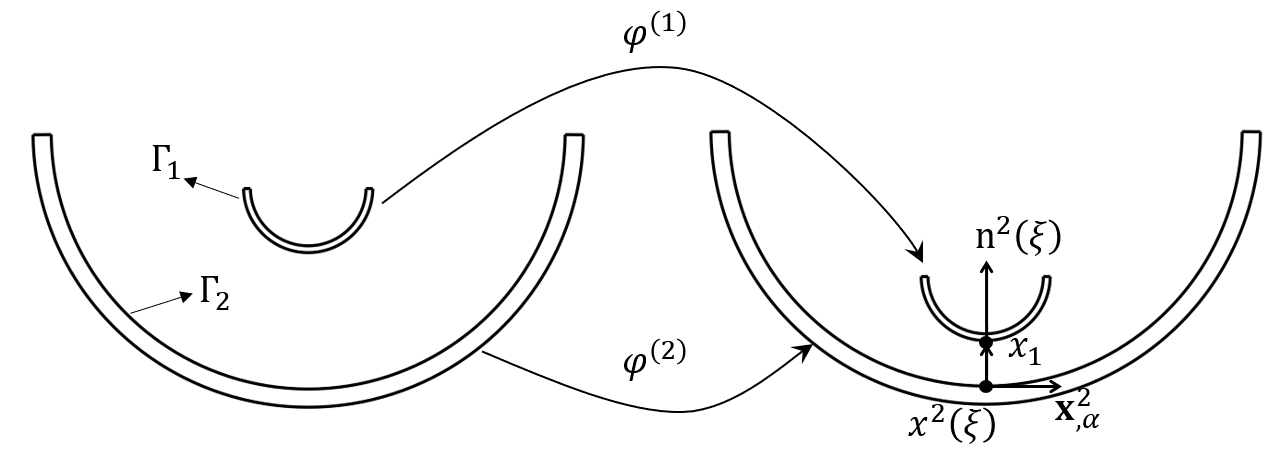}
\caption{The contact geometry}
\label{fig:contact geometry}
\end{figure}

The contact pairs mentioned above between the helical components and the cylinder layers alternatively come into stick and slip conditions with the cyclic bending loads. Therefore, the contact interface may be in normal contact or sliding contact depending upon the relative movement and friction properties. This present work considers the contact conditions without sliding. Figure \ref{fig:contact geometry} describes the contact geometries of outer sheath and steel tube. Both bodies are described in the initial configuration by $\Gamma_r$, where $r=1,2$ denotes contact bodies. $\mathbf{x}^{r}$ is the corresponding deformed configurations and $\varphi^r$ is the deformation mapping from the initial configuration to the deformed ones. Two contact bodies $\Gamma_r$ are described by the NURBS curves. Assume that the master body is $\mathbf{X}^{1}$ and slave body is $\mathbf{X}^{2}$. The minimum distance between two bodies is the most important parameter for contact detecting. Let each point in $\mathbf{X}^{1}$ to search the pairing point in $\mathbf{X}^{2}$ which is the most close to the master points to form contact pairs. Then the minimum distance between two points is: 
\begin{equation}
\|\mathbf{x}^{1}-\mathbf{x}^{2}\|=\min_{\mathbf{x}^{2}\subseteq\Gamma_{2}}\|\mathbf{x}^{1}-\mathbf{x}^{2}(\xi)\|    
\end{equation}

This minimizer is called the closest projection point of $\mathbf{x}^{1}$ and is denoted by $\mathbf{x}^{2}_p$. Then the relative position of two surfaces could be characterized by the gap. 

\begin{equation}
    g_N=[\mathbf{x}^1-\mathbf{x}^2(\xi)]\cdot\mathbf{n}^2(\xi)
    \label{eq:g_n}
\end{equation}

$\mathbf{n}^2(\xi)$ is the unit outward normal of $\mathbf{x}^{2}$ at $\mathbf{x}^{2}_p$. If $g_N>0$, the contact pairs separate. If $g_N=0$, two bodies come into contact. If $g_N<0$,  two bodies come into penetration.

\subsubsection{Contact detection}
Contact problems are well known nonlinear due to the inherent boundary constraints properties. The contact area between two bodies varies with the applied load. The contact conditions and relative positions of the bodies are not predetermined but instead evolve as the load varies. Therefore, the first challenge is how to identify the contact point and establish the contact pair. This issue could be resolved through trial-and-error iterations where contact point is iteratively determined. Once the contact point is identified, the corresponding contact constraints are imposed by incorporating the contact stiffness matrix. However, conventional node-to-surface (NTS) contact search methodologies employed in traditional finite element analyses exhibit limitations when applied to isogeometric analysis. This limitation arises from the fact that control points in NURBS geometric modeling do not align with interpolation nodes. To address this discrepancy, an extended approach known as knot-to-surface (KTS) contact search is commonly employed in isogeometric analysis. \\

In this work, the large steel tube is designated as the slave body, while the inner sheath layer serves as the master body. In the first load step, the contact is searched through each element within the slave body, evaluating each Gauss integration node vector $\mathbf{u}_{s}$. Subsequently, these node vectors are substituted into the NURBS curve Equation \ref{eq:NURBS} and this yields the corresponding node coordinates $\mathbf{u}_{s}(u_s)$. Then the potential contact coordinate $\mathbf{x}_{m}(u_m)$ is obtained based on the knot vector $u_m=0.5$. The tangential vector $\mathbf{x}_{\alpha}(u_m)$ is then determined by inserting into Equation 6. The closest point of the master body $\mathbf{x_m}(u_m)$ to the slave body could be identified when $err=0$. 

\begin{equation}err=[\mathbf{x_s}(u_s)-\mathbf{x_m}(u_m)]\cdot\mathbf{x_\alpha}(u_m)\end{equation}

But $err=0$ is extremely difficult to be obtained in the numerical analysis. Therefore, the $|err|<\epsilon $ is typically applied to approach the closest point $\mathbf{x_m}(u_m)$. $\epsilon $ could be set a positive and very small value, e.g. $1e^-10$. If $|err|>\epsilon$, an iteration should be conducted by the following process:

\begin{equation}u_\mathrm{m}=u_\mathrm{m}-\frac{err}{[\mathbf{x_s}(u_s)-\mathbf{x_m}(u_\mathrm{m})]\cdot\mathbf{x_\alpha}(u_\mathrm{m})}+\mathbf{x_\alpha}(u_\mathrm{m})\cdot\mathbf{x_\alpha}(u_\mathrm{m})\end{equation}

The iteration process is repeated until $|err|<\epsilon $. Then, the closest point on the master body $\mathbf{x_s}(u_s)$ can be identified and used to establish the contact pair. Subsequently, the normal vector $\mathbf{n^m}(u_m)$ at the coordinate  $\mathbf{x_s}(u_s)$ is obtained and substituted into the Equation \ref{eq:g_n} to determine the displacement function $g_n$ for this contact pair. This function is then utilized to ascertain whether the contact pair is activated.

\subsubsection{Penalty function}
The penalty method is widely recognized approach in addressing contact problems. The details of penalty method are described in nonlinear finite element methods\cite{kim_introduction_2014}. This section only explains the special aspects associated with the NURBS curve. The contact force is computed based on the Equation:

\begin{equation}P_{N}=\varepsilon<-g_{N}>
\label{eq:contactforce}
\end{equation}

The penetration is seen crucial to the contact force determination. However, when the contact boundaries are curves, the calculation error of penetration is unavoidable in the traditional finite element method. This is because the contact geometry can not be precisely described, leading to two types of errors as follows: \\

\begin{enumerate}
    \item In traditional finite element method, curves are approximately described by a number of straight lines. Figure \ref{fig:Penetrationduestraight} presents the difference between the actual penetration value and the approximation utilized in numerical analysis. The penetration in the numerical analysis is calculated based on the straight line connecting knots rather than accurately representing the curvature of the real curve. \\ 
\begin{figure}[htbp]
    \centering
    \includegraphics[width=0.6\linewidth]{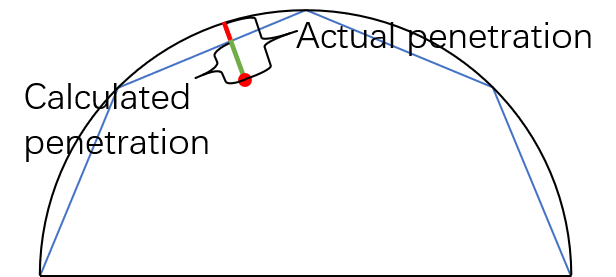}
    \caption{The penetration discrepancy between the actual value and the traditional finite element method}
    \label{fig:Penetrationduestraight}
\end{figure}
    \item The normal vector directions varies in the neighbouring elements since the approximated curve mentioned above is not continuous at the knot positions, as shown in Figure \ref{fig:Penetrationdue}. As a result, the penetrations obtained from adjacent elements vary at the same knot position. The corresponding contact forces at the same knot are therefore noncontinuous resulting in convergence problem.

\begin{figure}[htbp]
    \centering
    \includegraphics[width=0.6\linewidth]{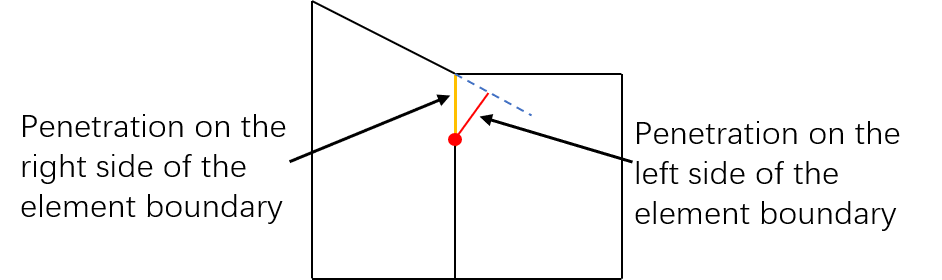}
    \caption{The normal vector directions varies in the neighbouring elements}
    \label{fig:Penetrationdue}
\end{figure}
\end{enumerate}

However, the contact geometry can be exactly described by the NURBS curve and then the penetration errors due to the geometry approximation can be eliminated. Moreover, the contact force remains continuous at knot positions, therefore, NURBS curve is capable of improving the calculation stability and efficiency. The isogeometric model offers significant advantages over traditional finite element methods in the contact analysis employing penalty function methods. 

\subsubsection{IGA contact algorithm scheme}
\begin{figure}[htbp]
    \centering
    \includegraphics[width=0.9\linewidth]{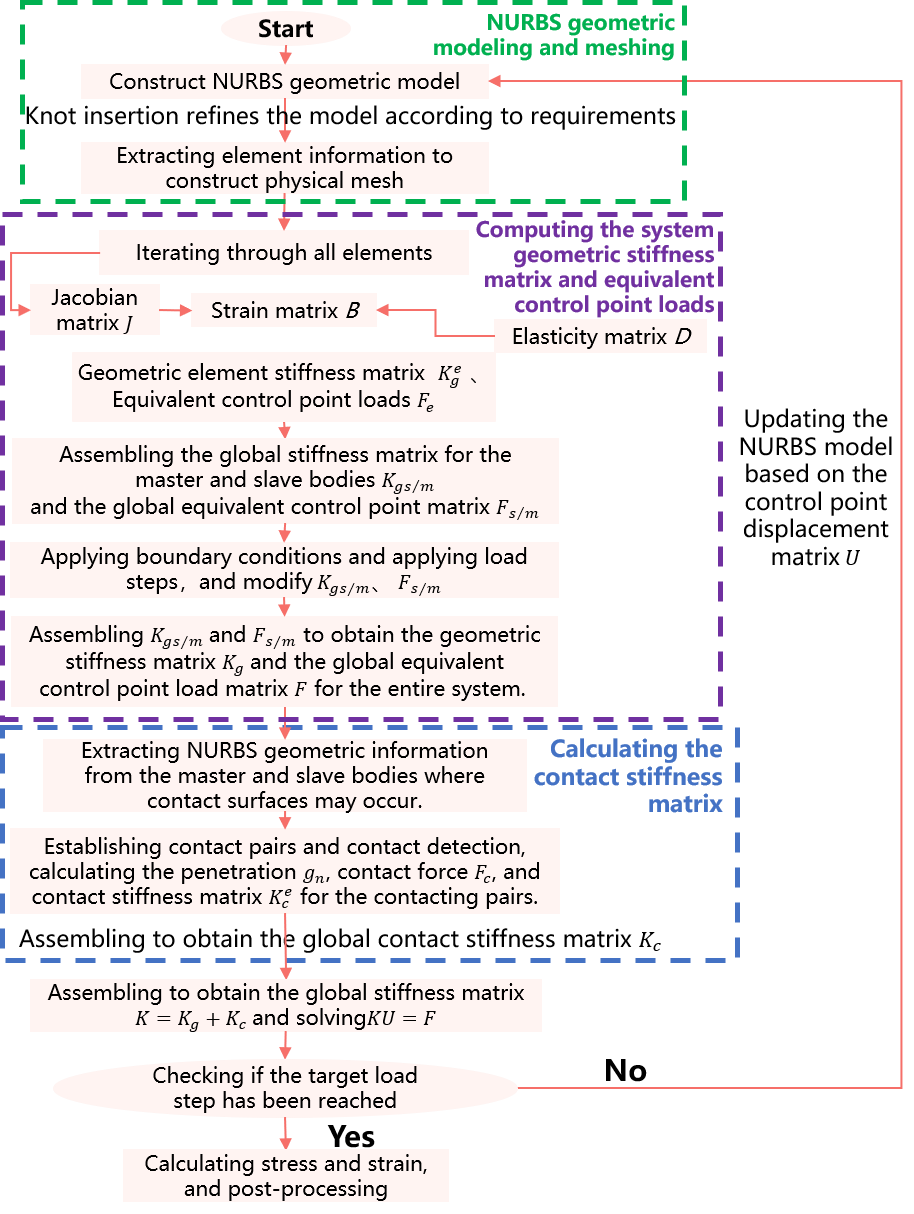}
    \caption{MATLAB program flow for isogeometric contact analysis}
    \label{fig:Program flow chart}
\end{figure}

The programming language employed in this study is MATLAB. The specific program flow is outlined as follows and shown in Figure \ref{fig:Program flow chart}:
\begin{enumerate}
    \item The NURBS geometric model of the specified object is established. This model includes the coordinates of the control points, the corresponding weight factors of the control points, the order of basis functions and the knot vectors. It will be updated based on the control point displacement matrix $\mathbf{U}$ when required in the specified load step.
    
    \item If the element refinement is required, inserting node is performed in the NURBS geometric model. This operation does not change the shape of the NURBS curve but only increases the number of elements. 
    
    \item The physical mesh is determined after defining the parameters of NURBS. Mesh generation is performed and involves extracting element information to construct a physical mesh conducive to subsequent calculations. 
    
    \item Iterations are performed in each element. In each iteration, the strain matrix $ \mathbf{B} $ and the Jacobian matrix $ \mathbf{J} $ are calculated for each element. Then, the element geometric stiffness matrix $ \mathbf{K_g^e} $ and the equivalent control point load matrix $ \mathbf{F^e} $ are obtained.
    
    \item The integration is performed within all elements with respect to the geometric stiffness matrices and equivalent control point load matrices.  This integration is carried out from the master and slave bodies into the global geometric stiffness matrix $ \mathbf{K_{gs/m}} $ and the global equivalent control point load matrix $ \mathbf{F_{s/m}} $.
    
    \item Boundary conditions are then applied. The global stiffness matrix $ \mathbf{K_{gs/m}} $ and the global equivalent control point load matrix $ \mathbf{F_{s/m}}$ are modified based on these boundary conditions requirements. They are modified as follows: the upper body are fixed in the X-axis direction to ensure that the slave body can only move along the y-axis direction for the slave body. Additionally, in each load step, a small increment of downward displacement is applied to the control points on the upper surface in each load step. For the master body, the lower surface elements are fixed in both x and y-axis directions to prevent sliding.
    
    \item Assembling the global geometric stiffness matrix $ \mathbf{K_g} $ and the global equivalent control point load matrix $ \mathbf{F} $ for the entire system by combining the global geometric stiffness matrix $ \mathbf{K_{gs/m}} $ and the global equivalent control point load matrix $ \mathbf{F_{s/m}} $ from the master and slave bodies.
    
    \item Extracting NURBS geometric information for the surfaces of the master and slave bodies where contact may occur.
    
    \item Iterating through every element and each Gaussian integration point on the slave body's contact surface, the nearest points on the master body's contact surface are found to form contact pairs. If the contact pairs are activated, the penetration $ g_n $, contact force $ F_c $, and contact stiffness matrix $ \mathbf{K_c^e} $ are computed. Based on the displacement $ d_m $ of nodes on the master body and the contact force, the contact stiffness $ S = F_c / d_m $ is calculated. It is crucial to note that the node displacement here refers to the displacement of coordinates on the physical mesh, not the displacement of control points. Lastly, the contact stiffness matrix $ \mathbf{K_c^e} $ is assembled into the overall contact stiffness matrix $ \mathbf{K_c} $.
    
    \item The global stiffness matrix $ \mathbf{K} $ is assembled by combining the global geometric stiffness matrix $ \mathbf{K_g} $ and the global contact stiffness matrix $ \mathbf{K_c} $. Subsequently, solving the Equation $ \mathbf{K}\mathbf{U} = \mathbf{F} $ yields the displacement matrix $ \mathbf{U} $ for all control points in the system.
    
    \item Repeat steps 1-10 until all load steps have been completed.
    
    \item Stress and strain of two deformed contact bodies can be obtained from this algorithm. More importantly, the contact surface properties including the normal contact stress, contact stiffness, penetrations can be extracted. 
\end{enumerate}

\section{Validation study}
In this section, the IGA contact algorithm is validated in terms of the numerical convergence capacity, efficiency, accuracy and stability. This validation was performed against the traditional finite element method in ABAQUS (named as ABAQUS contact model in the following text). Since this paper aims to investigate the inherent contact surface properties, the varying contact force along the helical length is not required to be considered. Therefore, the plane strain assumption is sufficient to establish the contact surface in both methods. The IGA contact algorithm has been introduced in Section 2. For the ABAQUS contact model, the CPE8R element is applied to mesh both the outer sheath and steel tube. This is an 8-node biquadratic element with the reduced integration to describe the plane strain deformation. \\

The contact pair between the outer sheath and large tube is chosen as an contact example for validations. The Young's modulus of large tube is $E_t=2×10^5 Mpa$, and Poison ratio is $\nu=0.3$. The external and internal diameters are 14.2 mm and 12.7 mm, respectively. The Young's modulus of outer sheath is $E_s=7×10^2 Mpa$, and Poison ratio is $\nu = 0.35$. The external and internal diameters are 60.45 mm and 58.45 mm, respectively. The penalty factor is selected to be 50 times of the Young's modulus of large tube, i.e., $\epsilon = 50E_t$. The vertical displacement is inserted on the surface of the large tube. The prescribed vertical displacement is $d=1\times 10^{-9} m$ in each load step and $D=1\time 10^{-7} m$ in total. \\

\subsection{Convergence capacity and efficiency}
Normally, the number of elements and DOFs determine the computation cost. Given that the IGA contact algorithm is implemented in MATLAB and differs from the ABAQUS solver, it is not appropriate to compare their computational efficiency based on runtime. Consequently, the element count and degrees of freedom (DOF) of two contact models are provided in Table \ref{tab:contactstresstable} for comparing their efficiency. In addition, Table \ref{tab:contactstresstable} shows the maximum normal stress comparison of two contact models. \\

\begin{figure}[htbp]
    \centering
    \includegraphics[width=\linewidth]{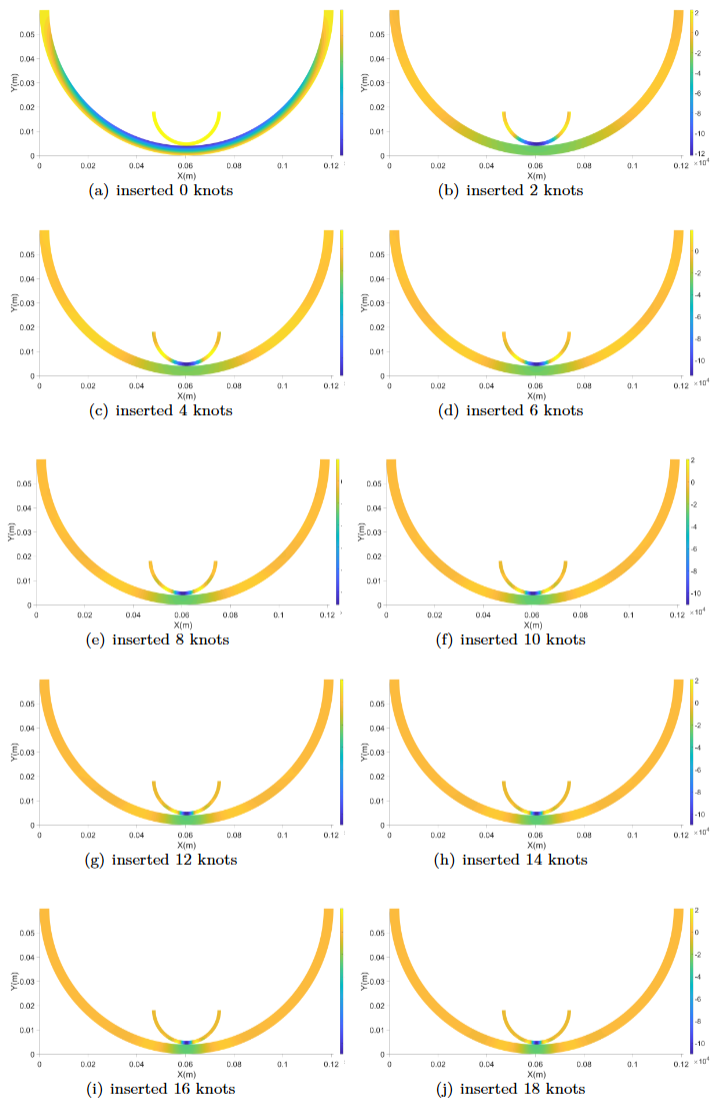}
    \caption{The normal stress distribution over the geometries of outer sheath and large tube by gradually inserting 18 knots}
    \label{fig:mesheffect}
\end{figure}

For the IGA contact algorithm, the original knot vector is $\mathbf U =\{0,0,0,1,1,1\}$ and the order of the basic functions is two in both x and y axis directions. 18 knots are gradually inserted into this original knot vector. Figure \ref{fig:mesheffect} shows the normal stress distribution of outer sheath and large tube by inserting different knots. The normal stress distribution around the contact area rapidly becomes intense after inserting four knots. The maximum normal stress is extremely efficiently to be converged when element number is 50 for the IGA contact algorithm. However, for the ABAQUS contact model, 44000 elements are required to obtain the same accurate normal stress and reasonable normal stress distribution, which demonstrates the excellent convergence capacity and efficiency of the IGA contact algorithm.

\subsection{Accuracy}
For the IGA contact algorithm, the maximum normal stress with 50 elements only has a 4\% deviation from the accurate value obtained when the element count is 882, which demonstrates the accuracy of IGA analysis. This is because the circular curve could accurately described even with few elements by using NURBS basic functions. Then the following validation study will further utilize this contact model with the insertion of 4 knots. \\

\begin{figure}[htbp]
    \centering
   \includegraphics[width=0.9\linewidth]{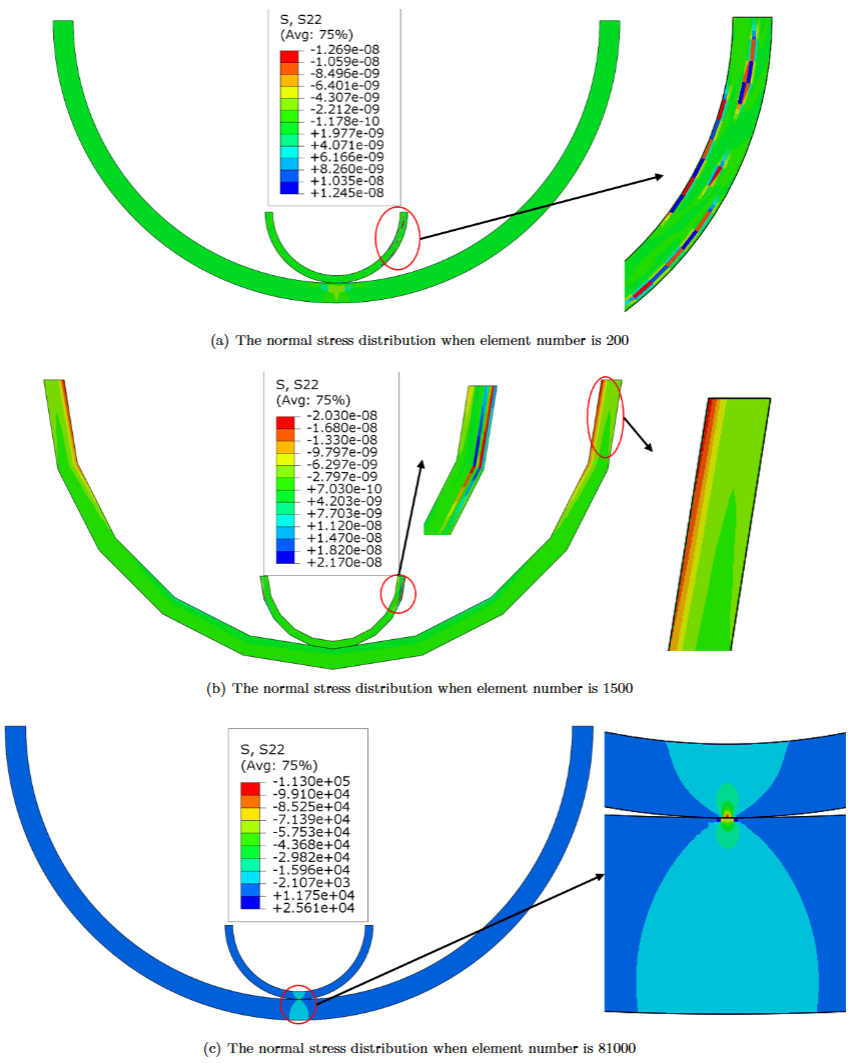}
    \caption{The normal stress distributions of large tube and outer sheath by using the ABAQUS contact model when element mesh is coarse and dense mesh}
    \label{fig:stresspropertyofABAQUS}
\end{figure}

Figure \ref{fig:stresspropertyofABAQUS} shows three typical normal stress distributions in the ABAQUS model for the element counts of 200, 1500 and 44000. With 200 elements, the normal stress is concentrated around the top area of the large tube. This area shows no contact between the large tube and outer sheath, which means an inaccurate stress prediction. As the element count increases, but remains below 52000, the normal stress distribution is chaos and occurs in the areas without contact. Furthermore, as seen in Table \ref{tab:contactstresstable}, the normal stress is nearly zero in the ABAQUS model when the element counts less than 52000.  This is because the contact detection and penetration are highly sensitive to the geometry. The traditional finite element method fail to accurately describe circular geometries with low element counts.  Once the mesh exceeds 52000 elements, the normal stress distribution aligns with that of the isogeometric model, though the convergence rate is slower. This demonstrates that the isogeometric method, with precise description of contact surfaces, is more efficient and accurate when applying the penalty function to calculate the contact force.  \\

\begin{table}[htbp]
    \centering
    \begin{tabular}{ p{1.8cm}<{\centering} p{1.4cm}<{\centering} p{1.8cm}<{\centering}  p{1.2cm}<{\centering} p{1.2cm}<{\centering} p{1.4cm}<{\centering}  p{1.8cm}<{\centering}  p{1.2cm}<{\centering} p{1.2cm}<{\centering}} 
        \hline
        \multicolumn{5}{c}{IGA contact algorithm}&\multicolumn{4}{c}{ABAQUS contact model}\\
        \hline
        Number of insert knot&Element numbers&Number of DOF      &$\sigma_{ymax}$(Pa)&Element numbers&Number of DOF &$\sigma_{ymax}$(Pa)\\
        0 &2  &36  &8.89e-6  &200    &1364     &1.12e7  \\
        2 &18 &100 &1.21e5   &1500   &9684    &2.26e-8 \\
        4 &50 &196 &1.15e5   &3000   &18764    &6.23e-8 \\
        6 &98 &324 &1.13e5   &6000   &37044   &7.16e-7 \\
        8 &162&484 &1.12e5   &13000  &79524   &2.12e-7 \\
        10&242&676 &1.10e5   &20000  &122644   &6.38e-8 \\
        14&450&1156&1.09e5   &26000  &158724   &1.00e-7 \\
        18&722&1764&1.10e5   &39000  &237924  &5.35e+4 \\
        20&882&2116&1.10e5   &44000  &268324  &1.06e+5 \\
        \hline
    \end{tabular}
    \caption{The maximum normal stress comparison of the large tube between the isogeometrical contact model and ABAQUS contact model when applying different element number}
    \label{tab:contactstresstable}
\end{table}

\subsection{Stability}

\begin{itemize}
\item Penalty factor effect
\end{itemize}

\begin{figure}[htbp]
\centering
\includegraphics[width=0.8\linewidth]{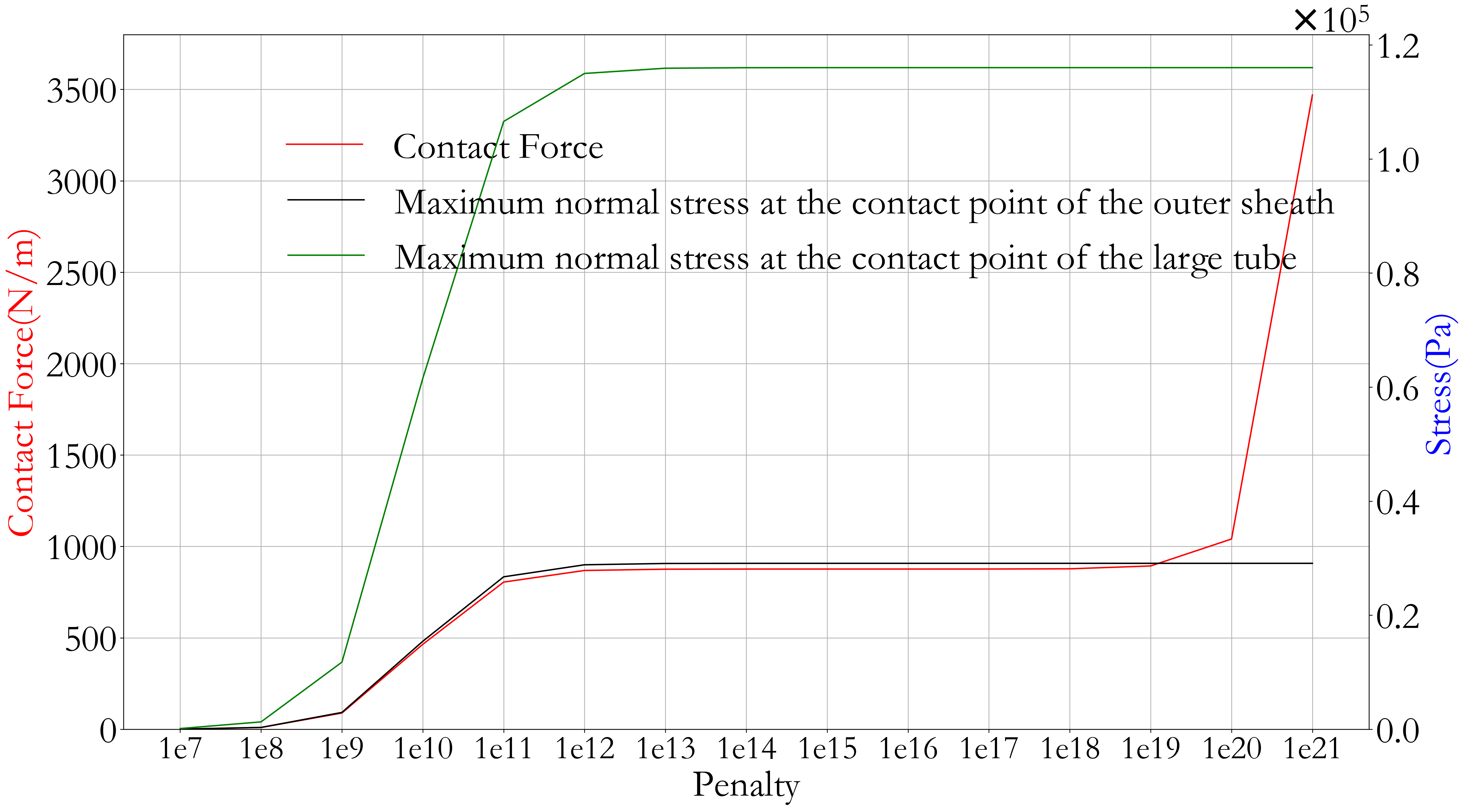}
\caption{The influence of penalty factor on contact force and maximum normal stress}
\label{fig:Penalty_force_stress} 
\end{figure}

\begin{figure}[htbp]
\centering
\includegraphics[width=0.8\linewidth]{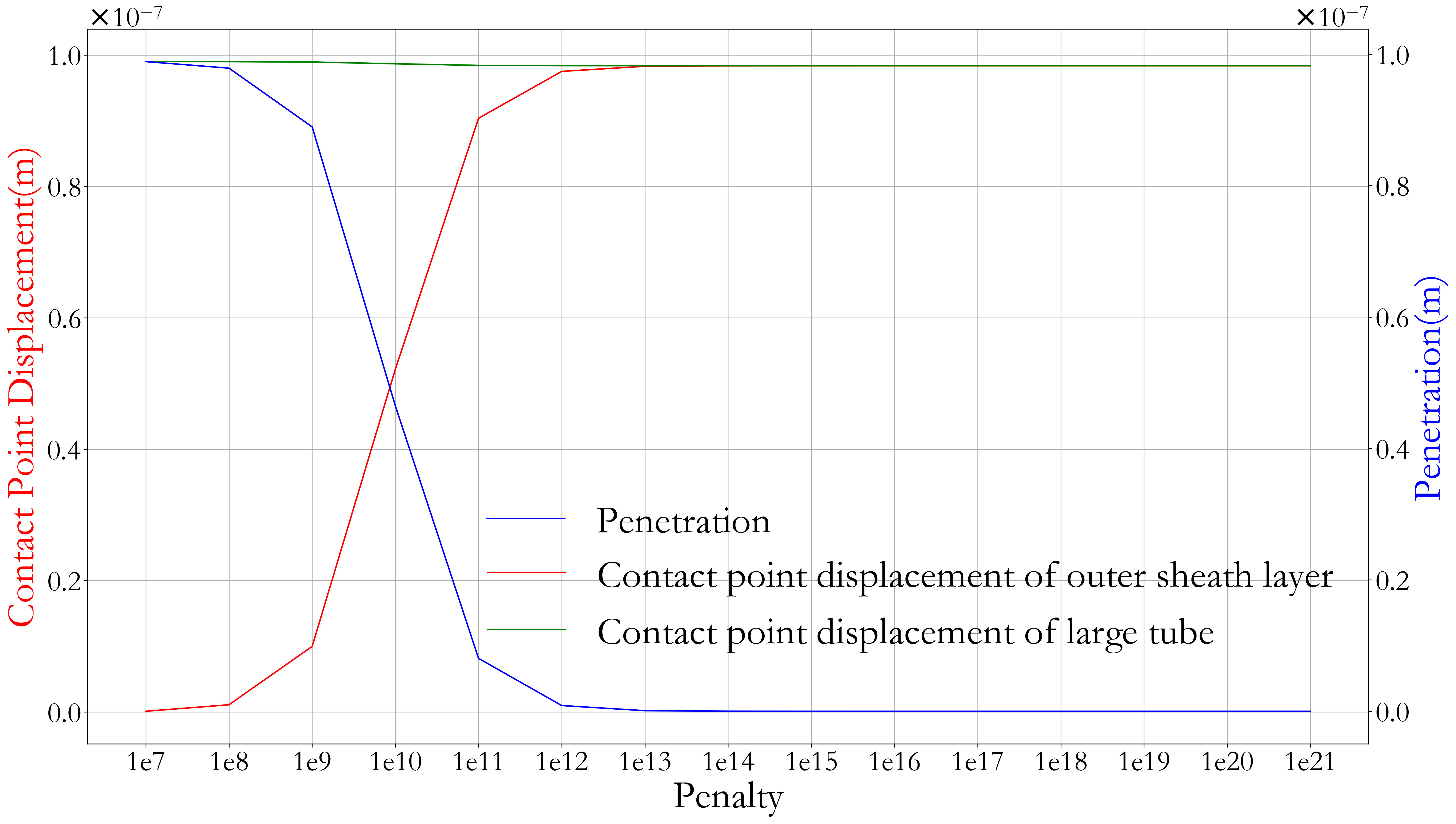}
\caption{The influence of penalty factor on contact point displacement and penetration}
\label{fig:Penalty_disp_penetra} 
\end{figure}

In contact problems employing penalty function methods, the penalty parameter should be sufficiently large to guarantee that the constraint is fulfilled. The value of the contact stiffness matrix is directly related to the magnitude of the penalty factor. When solving the Equation $\mathbf{K}\mathbf{U}=\mathbf{F}$, the inversion of the contact stiffness matrix may become highly unstable. Inversely, very small values within the inverted matrix can lead to matrix singularity, consequently affecting the stability and uniqueness of the solution. Therefore, in numerical simulations of contact problems, selecting an appropriate penalty factor is crucial. Typically, a sensitivity study is required to ensure that the penalty factor maintains numerical stability while accurately simulating the behavior of contact forces. Figure \ref{fig:Penalty_force_stress} and \ref{fig:Penalty_disp_penetra} illustrate the effect of the penalty factor on the contact behavior between a large steel tube and an outer sheath. The penalty factor of $1e12$ was found to be adequate in reflecting the constraint conditions, as the penetration gradually approaches zero. The contact displacements of the large tube and the outer sheath converge, indicating consistent displacement behavior. Furthermore, the maximum normal stress in both bodies stabilizes when the penalty factor exceeds $1e12$. 

\begin{itemize}
\item Relationship between contact force and displacement
\end{itemize}

\begin{figure}[htbp]
\centering
\includegraphics[width=0.5\linewidth]{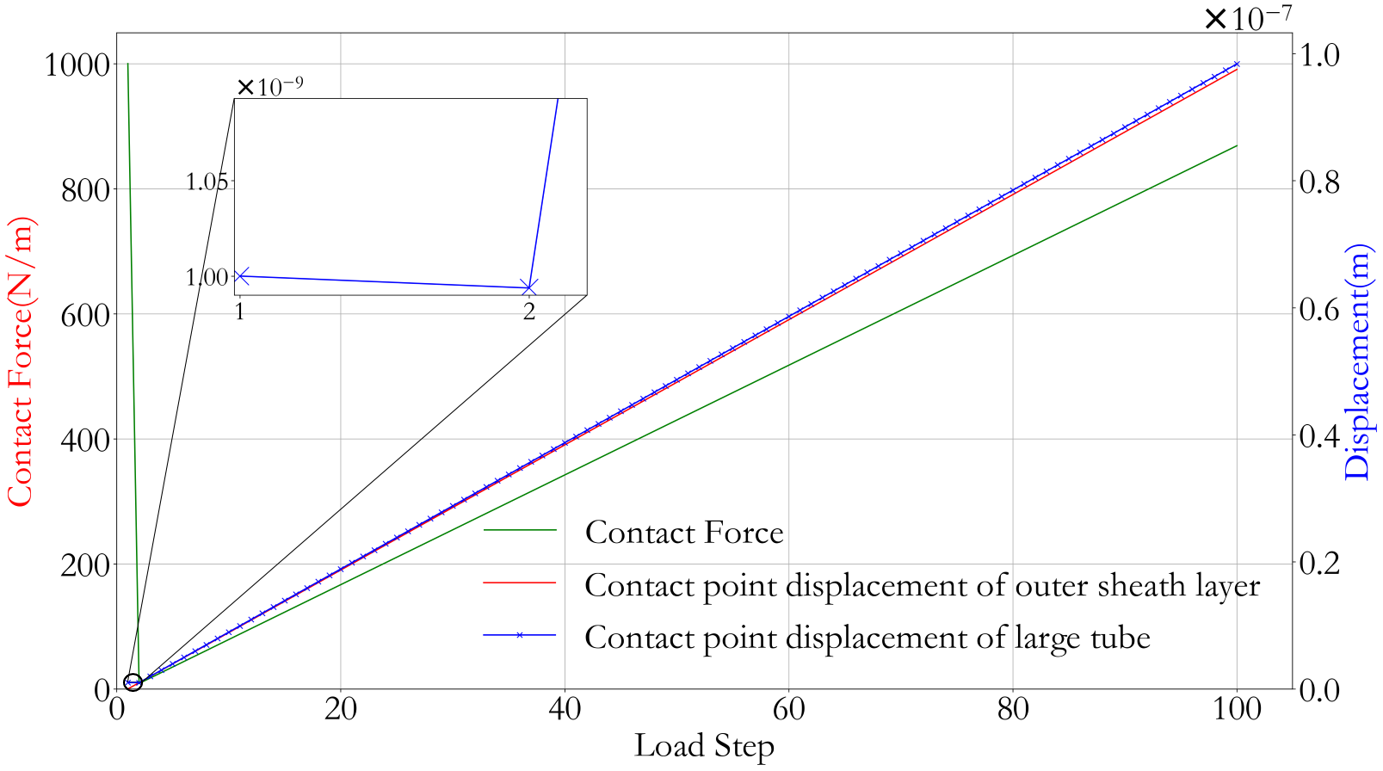}
\caption{The contact force and displacement of the contact surface between outer sheath and large tube }
\label{fig:contactvary} 
\end{figure}

Figure \ref{fig:contactvary} shows the contact point's displacement and contact force throughout the exposed vertical displacement history. In the first load step, the normal displacement is $1\times 10^{-9}m$ for the large tube and is zero for the outer sheath. Then the penetration is $1\times 10^{-9}m$ and the contact pair is activated. According to Equation \ref{eq:contactforce}, this yields a contact force is 1000$N/m$ as illustrated in this Figure \ref{fig:contactvary}. In the second load step, the contact stiffness matrix was calculated following the activation of the contact pair in the first load step. This stiffness matrix was then added to the geometry stiffness matrix to form the global stiffness matrix. Then the correct displacement of the control points at the contact area could be subsequently determined using the Equation $\mathbf{K}\mathbf{U}=\mathbf{F}$. The coordinates of the contact point were updated based on this newly obtained normal displacement. Then, the contact detection process was repeated and revealing that the penetration of the contact pair was significantly smaller than that in the first load step. Therefore, the contact force was found to be nearly zero and the normal displacement of large tube was reduced. More importantly, the contact force linearly increases with the increasing normal displacement of the contact point. \\

\begin{figure}[htbp]
    \centering
    \includegraphics[width=\linewidth]{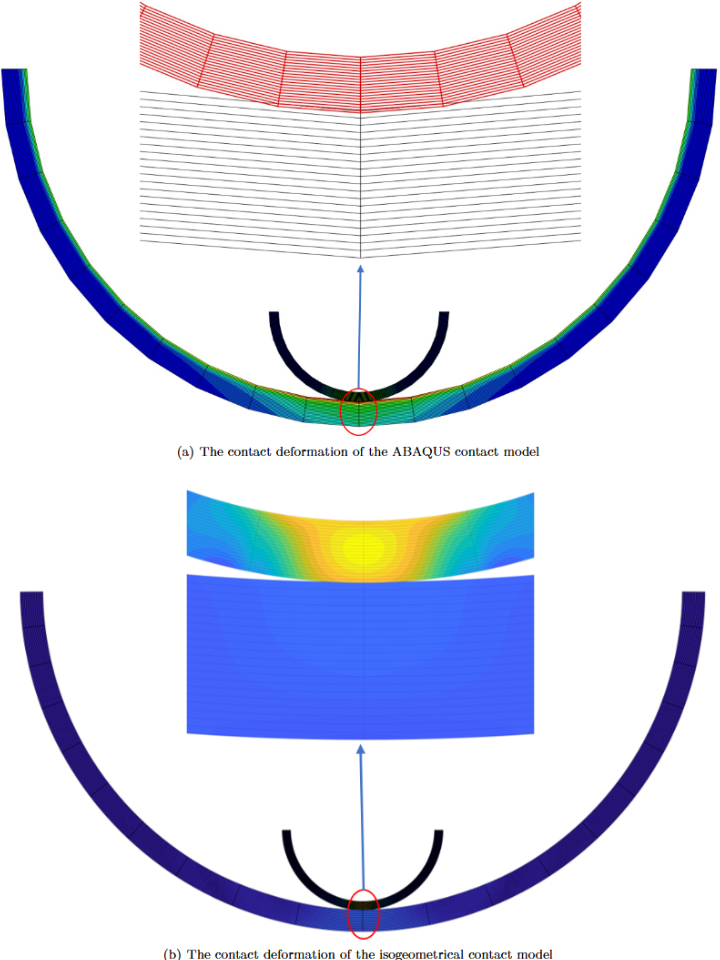}
    \caption{penetration comparison}
    \label{fig:penetrationcompare}
\end{figure}

\begin{itemize}
\item Contact area deformation
\end{itemize}
Figure \ref{fig:penetrationcompare} presents the contact area deformation of the ABAQUS contact model with the element count of 800 and the IGA contact algorithm with the element count of 50. For the ABAQUS contact model, the boundary at each knot between two neighbouring elements could be satisfied by the continuous condition $C^0$. However, the curvature at each knot is not continuous resulting in different directions of normal vectors and varying penetrations among neighbouring elements. Consequently, the contact detection would not be accurately described and results in chaotic stress distribution. In contrast, the contact area can be accurately described by the NURBS curve. This results in smaller penetrations than that in ABAQUS. The penetration is continuous and approaches the real condition, allowing the normal stress to converge rapidly to accurate values as seen in Figure \ref{fig:contactvary}. In addition, no vibration or noise is observed in IGA contact algorithm, indicating a highly stable calculation with respect to the convergence. \\

\section{Conclusions}
This paper presents the implementation of an Isogeometric Analysis (IGA) contact algorithm tailor-made to investigate the contact surface properties of dynamic umbilical and power cables. It requires the geometries and material properties of the contact pairs as inputs and delivers precise calculations of normal contact stress distribution, deformation, and contact stiffness within the contact area. The computed contact stiffness will be utilized in the subsequent stress analysis of dynamic umbilical and power cables, contributing to improved stability and accuracy in the numerical analysis, which forms the second part of this study.\\

The IGA contact algorithm has been validated against a finite element model established in ABAQUS, demonstrating good alignment in normal stress distribution and values between the two models. This validation highlights the algorithm's robust stability, convergence, and accuracy. Notably, the IGA contact algorithm achieves accurate results using only 50 elements, while the ABAQUS model requires over 44,000 elements to achieve similar accuracy, with only a 4\% deviation observed. The isogeometric method enables precisely modeling of the circular geometries of the contact pairs with a minimal element count by using NURBS curves. This approach achieves continuous curvature at knot positions, ensuring that penetration remains continuous, which enhances the accuracy of contact detection. As a result, the normal stress converges rapidly to accurate values.\\

\section*{Acknowledgement}
This research is supported by the National Natural Science Foundation of China (Grant No. 52001129) to the first author, Tianjiao Dai, and second author, Shuo Yang. The present work is also partially supported by the Major Project of Ningbo (Grant No.2021Z045) and the Fundamental Research Funds for the Central Universities HUST (Grant No.2019kfyXJJS006) to Tianjiao Dai.

\bibliographystyle{unsrt}
\bibliography{references_total}

\begin{thebibliography}{10}

\bibitem{li_efficient_2022}
Pengjie Li, Tianjiao Dai, Xing Jin, Leilei Dong, Shuaiqi Liu, Shuo Yang, Xianbo
  Xiang, and Hooi-Siang Kang.
\newblock An efficient fatigue analysis for the nonbondedflexible riser.
\newblock {\em Ships and Offshore Structures}, 17(10):2238--2253, October 2022.

\bibitem{li_finite_2021}
Xiaotian Li and Murilo~Augusto Vaz.
\newblock Finite element study on circular armour wire lateral buckling in
  umbilicals.
\newblock {\em Marine Structures}, 76:102895, March 2021.

\bibitem{li_finite_2023}
Xiaotian Li, Murilo~Augusto Vaz, and Anderson Barata~Custódio.
\newblock A finite element methodology for birdcaging analysis of flexible
  pipes with damaged outer layers.
\newblock {\em Marine Structures}, 89:103397, May 2023.

\bibitem{wang_full_2021}
Lidong Wang and Qianjin Yue.
\newblock A full layered numerical model for predicting hysteretic behavior of
  unbonded flexible pipes considering initial contact pressure.
\newblock {\em Applied Ocean Research}, 111:102626, 2021.
\newblock Publisher: Elsevier.

\bibitem{wang_novel_2022}
Lidong Wang, Naiquan Ye, and Qianjin Yue.
\newblock A novel helix contact model for predicting hysteretic behavior of
  unbonded flexible pipes.
\newblock {\em Ocean Engineering}, 264:112407, 2022.
\newblock Publisher: Elsevier.

\bibitem{zhang_research_2020}
Haichen Zhang, Lili Tong, Michael~Anim Addo, Jiaji Liang, and Lijun Wang.
\newblock Research on contact algorithm of unbonded flexible riser under
  axisymmetric load.
\newblock {\em International Journal of Pressure Vessels and Piping},
  188:104248, December 2020.

\bibitem{yin_experimental_2021}
Yuanchao Yin, Qingzhen Lu, Shanghua Wu, Zhixun Yang, Jun Yan, and Qianjin Yue.
\newblock Experimental study on the interlayer friction and wear mechanism
  between armor wires of umbilicals.
\newblock {\em Marine Structures}, 80:103102, 2021.
\newblock Publisher: Elsevier.

\bibitem{yang_study_2021}
Zhixun Yang, Qi~Su, Jun Yan, Shanghua Wu, Yandong Mao, Qingzhen Lu, and Hualin
  Wang.
\newblock Study on the nonlinear mechanical behaviour of an umbilical under
  combined loads of tension and torsion.
\newblock {\em Ocean Engineering}, 238:109742, October 2021.

\bibitem{qin_review_2024}
Xu~Qin, Mengmeng Zhang, Shixiao Fu, Huailiang Li, Jing Hou, and Yuwang Xu.
\newblock Review on researches and main influencing factors on mechanical
  properties of offshore wind power cables.
\newblock {\em Journal of Ocean Engineering and Science}, July 2024.

\bibitem{saevik_theoretical_2011}
Svein Sævik.
\newblock Theoretical and experimental studies of stresses in flexible pipes.
\newblock {\em Computers \& Structures}, 89(23-24):2273--2291, December 2011.

\bibitem{witz_axial-torsional_1992}
J.A. Witz and Z.~Tan.
\newblock On the axial-torsional structural behaviour of flexible pipes,
  umbilicals and marine cables.
\newblock {\em Marine Structures}, 5(2-3):205--227, January 1992.

\bibitem{saevik_stresses_1992}
Svein Sævik.
\newblock On stresses and fatigue in flexible pipes.
\newblock {\em PhD thesis}, Norwegian Institute of Technology, 1992.

\bibitem{dai_anisotropic_2018}
Tianjiao Dai, Svein Sævik, and Naiquan Ye.
\newblock An anisotropic friction model in non-bonded flexible risers.
\newblock {\em Marine Structures}, 59:423--443, 2018.
\newblock Publisher: Elsevier.

\bibitem{dai_friction_2017}
Tianjiao Dai, Svein Sævik, and Naiquan Ye.
\newblock Friction models for evaluating dynamic stresses in non-bonded
  flexible risers.
\newblock {\em Marine Structures}, 55:137--161, 2017.
\newblock Publisher: Elsevier.

\bibitem{dai_experimental_2020}
Tianjiao Dai, Svein Sævik, and Naiquan Ye.
\newblock Experimental and numerical studies on dynamic stress and curvature in
  steel tube umbilicals.
\newblock {\em Marine Structures}, 72:102724, 2020.
\newblock Publisher: Elsevier.

\bibitem{wu_hertzian_2016}
Chen-En Wu, Keng-Hui Lin, and Jia-Yang Juang.
\newblock Hertzian load–displacement relation holds for spherical indentation
  on soft elastic solids undergoing large deformations.
\newblock 97:71--76.

\bibitem{zhang_spherical_2014}
Man-Gong Zhang, Yan-Ping Cao, Guo-Yang Li, and Xi-Qiao Feng.
\newblock Spherical indentation method for determining the constitutive
  parameters of hyperelastic soft materials.
\newblock 13(1):1--11.

\bibitem{guo_modified_2020}
Zaoyang Guo, Meirong Hao, Li~Jiang, Dongfeng Li, Yang Chen, and Leiting Dong.
\newblock A modified hertz model for finite spherical indentation inspired by
  numerical simulations.
\newblock 83:104042.

\bibitem{el-abbasi_modelling_2001}
N.~El-Abbasi, S.~A. Meguid, and A.~Czekanski.
\newblock On the modelling of smooth contact surfaces using cubic splines.
\newblock 50(4):953--967.

\bibitem{de_lorenzis_mortar_2012}
L.~De~Lorenzis, P.~Wriggers, and G.~Zavarise.
\newblock A mortar formulation for 3d large deformation contact using
  {NURBS}-based isogeometric analysis and the augmented lagrangian method.
\newblock 49(1):1--20.

\bibitem{lu_circular_2009}
Jia Lu.
\newblock Circular element: {Isogeometric} elements of smooth boundary.
\newblock {\em Computer Methods in Applied Mechanics and Engineering},
  198(30):2391--2402, June 2009.

\bibitem{hegadekatte_multi_2010}
V.~Hegadekatte, J.~Hilgert, O.~Kraft, and N.~Huber.
\newblock Multi time scale simulations for wear prediction in micro-gears.
\newblock 268(1):316--324.

\bibitem{zhao_solution_2011}
Xin Zhao and Zili Li.
\newblock The solution of frictional wheel–rail rolling contact with a 3d
  transient finite element model: Validation and error analysis.
\newblock 271(1):444--452.

\bibitem{singh_analyses_2014}
Sarabjeet Singh, Uwe~G. Köpke, Carl~Q. Howard, and Dick Petersen.
\newblock Analyses of contact forces and vibration response for a defective
  rolling element bearing using an explicit dynamics finite element model.
\newblock 333(21):5356--5377.

\bibitem{xing_node--node_2019}
Weiwei Xing, Junqi Zhang, Chongmin Song, and Francis Tin-Loi.
\newblock A node-to-node scheme for three-dimensional contact problems using
  the scaled boundary finite element method.
\newblock 347:928--956.

\bibitem{otto_coupling_2019}
Peter Otto, Laura De~Lorenzis, and Jörg~F. Unger.
\newblock Coupling a {NURBS} contact interface with a higher order finite
  element discretization for contact problems using the mortar method.
\newblock 63(6):1203--1222.
\newblock Company: Springer Distributor: Springer Institution: Springer Label:
  Springer Number: 6 Publisher: Springer Berlin Heidelberg.

\bibitem{lu_isogeometric_2011}
Jia Lu.
\newblock Isogeometric contact analysis: Geometric basis and formulation for
  frictionless contact.
\newblock 200(5):726--741.

\bibitem{kim_introduction_2014}
Nam-Ho Kim.
\newblock {\em Introduction to Nonlinear Finite Element Analysis}.
\newblock Springer Science \& Business Media.
\newblock Google-Books-{ID}: {PDeSBQAAQBAJ}.

\end{thebibliography}
\end{document}